\documentclass[12pt]{amsart}
\usepackage{amscd, amsfonts, amssymb}
\usepackage{amscd}
\usepackage{epsfig}
\usepackage{amssymb}
\usepackage[mathscr]{eucal}
\usepackage{verbatim}
\usepackage{fullpage}
\usepackage{latexsym}
\usepackage{lscape}

\newcommand\cc{\mathfrak c}
\renewcommand\gg{\mathfrak g}
\newcommand\hh{\mathfrak h}

\newcommand\mm{\mathfrak m}

\newcommand\inverse{{^{-1}}}
\renewcommand{\check}{^{\vee}}

\newcommand\DD{\mathcal D}
\newcommand\NN{\mathcal N}
\newcommand\UU{\mathcal U}

\newcommand\ZZ{\mathbb Z}

\DeclareMathOperator{\Ad}{Ad}

\DeclareMathOperator{\Char}{char}

\DeclareMathOperator{\GL}{GL}
\DeclareMathOperator{\SL}{SL}

\DeclareMathOperator{\SO}{SO}
\DeclareMathOperator{\SP}{Sp}
\DeclareMathOperator{\Hom}{Hom}
\DeclareMathOperator{\Lie}{Lie}

\DeclareMathOperator{\rank}{rank}

\numberwithin{equation}{section}

\newtheorem{thm}[equation]{Theorem}

\newtheorem{lem}[equation]{Lemma}
\newtheorem{cor}[equation]{Corollary}
\newtheorem{prop}[equation]{Proposition}

\theoremstyle{definition}
\newtheorem{defn}[equation]{Definition}
\newtheorem{exmp}[equation]{Example}    
\theoremstyle{remark}
\newtheorem{rem}[equation]{Remark}
    
\theoremstyle{remark}
\newtheorem{rems}[equation]{Remarks}

\thanks{2000 {\it Mathematics Subject Classification}. 
20G15, 14L30 (17B50)} 
\keywords{Cocharacters associated to nilpotent elements}

\title[cocharacters associated to nilpotent elements]
{On cocharacters associated to nilpotent elements of reductive 
groups}

\author[R.\  Fowler]{Russell Fowler} 
\address
{School of Mathematics, University of Birmingham, 
Birmingham, B15 2TT, UK}  

\author[G.\ R\"ohrle]{Gerhard R\"ohrle}
\address
{Fakult\"at f\"ur Mathematik,
Ruhr-Universit\"at Bochum,
Universit\"atsstrasse 150,
D-44780 Bochum, Germany}
\email{gerhard.roehrle@rub.de}

\begin{document}

\begin{abstract}
Let $G$ be a connected reductive linear algebraic group
defined over an algebraically closed field of characteristic $p$.
Assume that $p$ is good for $G$.
In this note we consider particular classes of 
connected reductive subgroups $H$ of $G$ 
and show that the cocharacters of $H$ that are associated to a given
nilpotent element $e$ in the Lie algebra of $H$ 
are precisely the cocharacters of $G$ associated to $e$ that take
values in $H$.
In particular, we show that this is the case provided $H$
is a connected reductive subgroup of $G$ of maximal rank;
this answers a question posed by J.C.~Jantzen.
\end{abstract}

\maketitle


\section{Introduction}
\label{s:intro}
Let $G$ be a connected reductive linear algebraic group
defined over some algebraically closed field $k$,
let $\gg$ be its Lie algebra and $\NN$ be the nilpotent cone of $\gg$. 
The Jacobson--Morozov Theorem allows one to associate an 
$\mathfrak{sl}(2)$-triple to any given non-zero nilpotent element in $\NN$
in characteristic zero or large positive characteristic.
This is an indispensable tool in the 
Dynkin--Kostant classification of the
nilpotent orbits in characteristic zero as well as 
in the Bala--Carter classification of unipotent conjugacy classes of $G$
in large prime characteristic, see \cite[\S 5.9]{carter}.
In good characteristic there is a replacement for
$\mathfrak{sl}(2)$-triples, so called \emph{associated 
cocharacters}; see Definition \ref{def:assoc} below. 
These cocharacters have become important tools in the classification theory of 
unipotent and nilpotent classes of reductive algebraic groups
in good characteristic, see
for instance \cite[\S 5]{Jantzen}, \cite{mcninch}, \cite{mcninchsommers},
\cite{premet2},  and \cite{premet} for more details.

In \cite[\S 5.6]{Jantzen}, J.C.~Jantzen studies 
the behaviour of 
cocharacters associated to nilpotent elements under 
elementary operations of algebraic groups such 
as passing to derived subgroups, taking direct products and isogenies
of reductive groups.
While these cocharacters 
also behave well with respect to inclusions of reductive subgroups
in characteristic zero, this is not the case in general 
in positive characteristic, \cite[\S 5.12]{Jantzen};
given that they serve as a replacement for $\mathfrak {sl}(2)$-triples, 
this is somewhat surprising. 

More precisely, in \cite[Claim 5.12]{Jantzen}, 
Jantzen showed that if $\Char k = 0$
and $H$ is a connected reductive subgroup of $G$ with 
Lie algebra $\hh \subseteq \gg$, then 
\begin{itemize}
\label{dagger}
\item[($\dagger$)]
the cocharacters of $H$ associated to $e \in \hh \cap \NN$ 
are precisely  
the cocharacters of $G$ associated to $e$ that take values in $H$.
\end{itemize}
Jantzen continues to give an 
example in \cite[\S 5.12]{Jantzen} which shows
that ($\dagger$) fails in general in positive characteristic, 
even when $\Char k$ is good for $G$.

G.\ McNinch pointed out that this failure ultimately stems from the fact
that the representations of $H$ are not semisimple in general
in positive characteristic, as is the case in Jantzen's counterexample.
The following construction due to McNinch shows
that $(\dagger)$ fails generically: 
Let $H$ be a connected reductive group whose Coxeter number satisfies 
$h < p = \Char k$.
Let $K \le H$ be a \emph{principal} $\SL_2$-subgroup of $H$ 
and let $V$ be a faithful representation of $H$ for which 
there is a $K$-composition factor of the restriction of $V$ to $K$ having
a non-restricted highest weight.
See \cite{mcninch2} for a definition, existence and uniqueness
up to conjugacy of principal $\SL_2$-subgroups, see also \cite{seitz2}.
A maximal torus of $K$ determines a cocharacter $\lambda$ of $H$ associated 
to some 
nilpotent element $e$ of $\Lie K$. 
However, viewed as a cocharacter of $\GL(V)$, 
$\lambda$ is not associated to $e$.
In this sense, $(\dagger)$ fails for $H \le \GL(V)$
for most representations $V$ of $H$.

Nevertheless, a calculation of Jantzen shows that for $G$ of classical type
and $\Char k$ a good prime for $G$,
$(\dagger)$ does hold provided the subgroup $H$ of $G$ is of maximal rank,
see \cite[\S 5.12]{Jantzen}. 
In this note we give a general proof showing that this is indeed
always the case irrespective of the type of $G$:

\begin{thm}
\label{thm:maxrank0}
Suppose $\Char k$ is good for $G$.
Let $H$ be a connected reductive subgroup of $G$ of maximal rank.
Then $(\dagger)$ holds.
\end{thm}

This answers a question posed by J.C.~Jantzen, see \cite[\S 5.12]{Jantzen}.
We then extend this result 
to the case when $H$ is a regular reductive subgroup of $G$,
see Theorem \ref{thm:regular}.

In Theorem \ref{thm:reduction} we show that $(\dagger)$ holds provided
there is at least one cocharacter of $G$ that is associated to $e$ 
and lies in $H$. In our principal result we show that  
this is the case for a special class of subgroups:

\begin{thm}
\label{thm:main0}
Let $G$ be a connected reductive algebraic group. 
Suppose that $\Char k$ is a good prime for $G$.
Let $S$ be a linearly
reductive group acting on $G$ by automorphisms and 
set $H = C_G(S)^\circ$. 
Then $(\dagger)$ holds.
\end{thm}

Theorem \ref{thm:main0} is a consequence of 
Theorem \ref{thm:reduction} and Lemma \ref{l:S-fixedtau}.
Our proof uses ideas due to R.W.~Richardson 
that he used in order to show that 
any $S$-stable parabolic subgroup of $G$ admits an 
$S$-stable Levi subgroup in \cite[Prop.\ 6.1, \S 6.2, \S 6.3]{rich0}.

In turn Theorem \ref{thm:maxrank0} is a direct consequence of 
Theorem \ref{thm:main0}. 

In Subsection \ref{sub:local} we consider arbitrary connected 
reductive subgroups $H$ of $G$ and give criteria for 
a given nilpotent element $e\in \hh$ that will ensure $(\dagger)$ to hold.
For instance, in Lemma \ref{lem:dist} we prove that this is the case 
provided $e\in \hh$ is distinguished in $\gg$.
In another result,
Theorem \ref{thm:localrank}, we show that $(\dagger)$ holds
given that the ranks of the centralizers of $e$ in $H$ and $G$ coincide.

We assume throughout that $\Char k$ is good for $G$. 
Then there exists a $G$-equivariant homeomorphism $\NN \to \UU$ 
between the nilpotent cone $\NN$ of $\gg$ and the unipotent variety $\UU$ 
of $G$. Such a map is called a \emph{Springer isomorphism},
see \cite[III, 3.12]{SS} and \cite[Cor.\ 9.3.4]{BaRi}.
By means of such a map, all the results below admit analogues for 
associated cocharacters of unipotent elements in $G$.

For general results on algebraic groups we cite Borel's book \cite{borel},
and for basic results on cocharacters associated to nilpotent elements, 
we refer the reader to Jantzen's monograph \cite[Ch.\ 5]{Jantzen}
and the articles by McNinch--Sommers \cite{mcninchsommers}
and Premet \cite{premet}.
\goodbreak

\vfil
\eject

\section{Preliminaries} 
\label{s:prelim}

\subsection{Notation} 
\label{sub:not}
Throughout, $G$ is a connected reductive algebraic group defined over 
an algebraically closed field $k$ and $p = \Char k$ is a good prime for $G$,
although many results hold without this assumption.
We denote the Lie algebra of $G$ by $\Lie G$ or by $\gg$; 
likewise for closed subgroups of $G$.
For $e \in \gg$ and $g \in G$ we denote the adjoint action of 
$g$ on $e$ by $\Ad(g)e$.
The centralizers of $e$ in $G$ and $\gg$ are
$C_G(e) = \{g\in G \mid \Ad(g) e = e\}$
and $\cc_\gg(e) = \{x \in \gg \mid [x,e] = 0\}$, respectively.
We write $Z(G)$ for the centre of $G$.

Let $H$ be a closed subgroup of $G$. We write $H^{\circ}$ for the 
identity component of $H$ and 
$C_G(H) = \{ g \in G \mid ghg\inverse = h \text{ for all } h \in H\}$
for the centralizer of $H$ in $G$. 
The normalizer of $H$ in $G$ is $N_G(H) = \{g \in G \mid gHg\inverse = H\}$.
The derived subgroup of $H$ is denoted by $\DD H$ and we write
$\rank H$ for the dimension of a maximal torus of $H$.
The unipotent radical of $H$ is denoted by $R_u(H)$.
A \emph{Levi subgroup} of $H$ is a complement to $R_u(H)$ in $H$, 
\cite[Defn.\ 11.22]{borel};
in contrast to \emph{loc.\ cit.}, we do not require $H$ to be connected;
we also refer to the semi-direct product 
of a Levi subgroup of $H$ and $R_u(H)$ as 
a \emph{Levi decomposition} of $H$.

By a Levi subgroup of $G$ we mean a Levi subgroup of a 
parabolic subgroup of $G$.
Let $S$ be a torus of $G$. Then $C_G(S)$ is a Levi subgroup of $G$,
\cite[Thm.\  20.4]{borel}. 
Note that $C_G(S)$ is connected, \cite[Cor.\ 11.12]{borel}.
Moreover, by \cite[Prop.\  8.18]{borel} there exists an element
$s \in S$ such that $C_G(s) = C_G(S)$. 
Conversely, every Levi subgroup of $G$ is of this form, e.g., see
Lemma \ref{lem:cochars}(ii).

Let $Y(G) = \Hom(k^*,G)$ denote the set of cocharacters of $G$. 
For $\mu \in Y(G)$ we write $C_G(\mu)$ for $C_G(\mu(k^*))$.
For $\mu\in Y(G)$ 
and $g\in G$ we define the \emph{conjugate cocharacter}
$g\cdot \mu \in Y(G)$ by 
$(g\cdot \mu)(t) = g\mu(t)g\inverse$; 
this gives a left action of $G$ on $Y(G)$.
For $H$ a (connected) reductive subgroup of $G$, let $Y(H) = \Hom(k^*,H)$
denote the set of cocharacters of $H$. There is an obvious 
inclusion $Y(H) \subseteq Y(G)$.

If $S$ is a linear algebraic group acting on $G$ by automorphisms, 
then we say that $G$ is an \emph{$S$-group}; we write $s \cdot g$ 
for the action of $s \in S$ on $g \in G$.
The subgroup of $G$ consisting of the $S$-fixed points is denoted by
$C_G(S) = \{g \in G \mid s\cdot g = g  \text{ for all } s \in S\}$.
If $G$ is an $S$-group, then $S$ acts naturally by means of 
Lie algebra automorphisms on $\gg$.
By abuse of notation, we simply denote the action of $S$ on $\gg$
by $s \cdot e$ for $s \in S$ and $e \in \gg$.
We denote the subalgebra of $\gg$ consisting 
of the $S$-fixed points for the induced action on $\gg$ by  
$\cc_\gg(S) = \{e \in \gg \mid s\cdot e = e  \text{ for all } s \in S\}$.
Also, $S$ acts on $Y(G)$ by acting on the image of a cocharacter in $G$:
$(s \cdot\lambda)(t) = s \cdot \lambda(t)$ for $s \in S$, $\lambda \in Y(G)$,
and $t \in k^*$.

More generally, if $S$ and the $S$-group $G$ both act morphically on 
an algebraic variety $X$, then, following \cite[(2.1)]{rich0},
the actions of $G$ and $S$ are said to be \emph{compatible}, provided
\begin{equation}
\label{e:compatible}
s \cdot (g\cdot x) = (s\cdot g) \cdot (s \cdot x)
\end{equation}
for all $s \in S$, $g \in G$ and $x \in X$.
This is the unique action that defines a morphic action of the
semidirect product of $G$ and $S$ on $X$ which extends the actions of 
both $G$ and $S$ on $X$, see \cite[\S 2]{rich0}.
All actions by an $S$-group $G$ together with $S$ 
considered in this paper are compatible in the sense
of \eqref{e:compatible}.

Let $T$ be a maximal torus of $G$.
Let $\Psi = \Psi(G,T)$ denote the set of roots of $G$
with respect to $T$.
Fix a Borel subgroup $B$ of $G$ containing $T$ and let  
$\Pi = \Pi(G, T)$   
be the set of simple roots of $\Psi$ defined by $B$. Then
$\Psi^+ = \Psi(B)$ is the set of positive roots of $G$.
For $\beta \in \Psi^+$ write 
$\beta = \sum_{\alpha \in \Pi} c_{\alpha\beta} \alpha$
with $c_{\alpha\beta} \in \mathbb N_0$.
A prime $p$ is said to be \emph{good} for $G$
if it does not divide $c_{\alpha\beta}$ for any $\alpha$ and $\beta$, 
\cite[Defn.\ 4.1]{SS}.

\subsection{Linearly reductive groups}
\label{sub:linred}
We refer to Richardson's article \cite{rich0} for information on 
centralizers of the action of linearly reductive groups on connected 
reductive groups.
Recall that a linear algebraic group $S$, not necessarily connected, 
is said to be {\em linearly reductive} if every rational representation 
of $S$ is semisimple.  It is well known that in characteristic zero, 
$S$ is linearly reductive if and only if $S^\circ$ is reductive.  
In characteristic $p > 0$, $S$ is linearly reductive if and only 
if every element of $S$ is semisimple if and only if $S^\circ$ 
is a torus and $|S/S^\circ|$ 
is coprime to $p$, see \cite[\S4, Thm.\  2]{nagata}.
In particular, a torus is linearly reductive.

In the sequel we require some fundamental results concerning centralizers 
of linearly reductive groups acting on connected reductive groups; 
the following facts are \cite[Lem.\ 4.1, 
Prop.\ 10.1.5]{rich0}.

\begin{prop}
\label{prop:linear1}
Let $G$ be a connected reductive algebraic group 
and $S$ a linearly reductive algebraic 
group acting on $G$ so that $G$ is an $S$-group. Then we have
\begin{itemize}
\item[(i)]
$C_G(S)^\circ$ is reductive;
\item[(ii)]
$\Lie C_G(S) = \cc_\gg(S)$.
\end{itemize}
\end{prop}

The following result is due to R.W.\ Richardson, cf.\ \cite[Prop.\ 6.1]{rich0}.

\begin{prop}
\label{prop:richardson}
Let $G$ be a connected reductive algebraic group 
and $S$ a linearly reductive algebraic 
group acting on $G$ so that $G$ is an $S$-group. 
Let $K$ be a (not necessarily connected)
closed $S$-stable subgroup of $G$.
Suppose that $K$ admits a Levi decomposition 
such that $R_u(K)$ acts simply transitively 
on the set of all Levi subgroups of $K$.
Then $K$ admits an $S$-stable Levi subgroup.
\end{prop}

\begin{proof}
Although \cite[Prop.\ 6.1]{rich0} is  
only stated for parabolic subgroups of $G$, Richardson's proof applies in this
slightly more general setting \emph{mutatis mutandis}.
The given conditions on the Levi subgroups of $K$ are precisely 
the relevant properties of the set of Levi subgroups 
of a parabolic subgroup of $G$ 
in Richardson's proof, see \cite[\S 6.3]{rich0}.
In addition, a general result on the vanishing of 
the (non-commutative) cohomology group  
$H^1(S,R_u(K))$ is needed; this is
proved by Richardson in \cite[Lem.\ 6.2.6]{rich0}.
\end{proof}

\begin{rem}
\label{rem:jantzen11.24}
If $S$ is a subgroup of $K$ in Proposition \ref{prop:richardson}, 
then the conclusion 
is that there exists a Levi subgroup 
of $K$ containing $S$,
see Jantzen's generalization \cite[Lem.\ 11.24]{Jantzen}
of Mostow's theorem \cite[Thm.\ 7.1]{mostow}.
In that case the assumption that  
$R_u(K)$ acts simply transitively on the set of 
all Levi subgroups of $K$ is not necessary.

The essence of the arguments in the proofs of 
\cite[Lem.\ 11.24]{Jantzen},
\cite[Thm.\ 7.1]{mostow}, and 
\cite[Lem.\ 6.2.6]{rich0} is the 
vanishing of the cohomology $H^1(F,U)$, where 
$F$ is a finite group whose order is coprime to $p$ and 
$U$ is a (finite) unipotent group.
\end{rem}

\begin{rem}
\label{rem:ssauto}
Let $\sigma$ be a semisimple automorphism
of $G$, \cite[\S 7]{St}. Then there is an embedding 
$G \le \GL_n$ of algebraic groups for some $n$ such that $\sigma$ is given 
by conjugation by a semisimple element, say $s$, in $\GL_n$.
Thus $s$ belongs to some maximal torus of $\GL_n$, so 
the algebraic subgroup $S$ of $\GL_n$ generated by $s$ 
consists of semisimple elements. Thus
$S$ is linearly reductive, by \cite[\S 4, Thm.\ 2]{nagata}.
Clearly, $S$ depends on the choice of the ambient group $\GL_n$, but the 
fact that $S$ is linearly reductive does not.
\end{rem}

\subsection{Regular reductive subgroups}
\label{sub:regular}

Let $H$ be a closed (not necessarily connected)  
subgroup of $G$ normalized by some maximal torus $T$
of $G$; that is, 
a \emph{regular} subgroup of $G$ (reductive 
regular subgroups are often also referred to as 
\emph{subsystem subgroups} in the literature).
In this case the root spaces of $\mathfrak h$ relative to $T$ 
are also root spaces of  $\mathfrak g$ relative to $T$, 
and the set of roots of $H$ with respect to $T$, 
$\Psi(H) = \Psi(H, T) =  \{\alpha \in \Psi \mid \mathfrak g_\alpha
\subseteq \mathfrak h\}$,  is a subset of $\Psi$, 
where $\mathfrak g_\alpha$ denotes the root space in $\mathfrak g$ 
corresponding to $\alpha$.
If $\Char k$ does not divide any of the structure constants
of the Chevalley commutator relations of $G$,
then $\Psi(H)$ is closed under addition. 
In particular, this is the case when $\Char k$ is a good prime for $G$. 
If $H$ is reductive and regular, then 
$\Psi(H)$ is a semisimple subsystem of $\Psi$. 

Recall that for $s \in G$ semisimple, $H = C_G(s)^{\circ}$ 
is called a \emph{pseudo-Levi subgroup} of $G$, 
cf.\ \cite[\S 6]{mcninchsommers}. 
Since $s$ is contained in a maximal torus $T$ of $G$, it follows that $H$ 
is regular of maximal rank. Moreover, 
$\Psi(H) = \Psi(H, T) =  \{\alpha \in \Psi(G,T) \mid \alpha(s) = 0\}$.

\subsection{Kempf--Rousseau Theory}
\label{sub:git}

Let $G$ be a reductive group acting on an affine variety $X$. 
For $x \in X$ let $G \cdot x$ denote the $G$-orbit of $x$ in $X$ and 
$C_G(x)$ the stabilizer of $x$ in $G$.  
Let $\phi : k^* \to X$ be a morphism of algebraic varieties. 
We say that  $\underset{t\to 0}{\lim}\, \phi(t)$ exists 
if there exists a morphism $\widehat\phi :k \to X$ 
(necessarily unique) whose restriction to $k^*$ 
is $\phi$; if this limit exists, then we set
$\underset{t\to 0}{\lim}\, \phi(t) = \widehat\phi(0)$.

Recall the characterization of parabolic subgroups of $G$
in terms of cocharacters of $G$, e.g.\ see \cite[Prop.\  8.4.5]{spr2}.

\begin{lem} 
\label{lem:cochars} 
Given a parabolic subgroup $P$ of $G$ and any Levi subgroup $L$ of $P$,
there exists $\lambda \in Y(G)$ 
such that the following hold:
\begin{itemize}
\item[(i)]   $P = P_\lambda := \{g\in G \mid \underset{t\to 0}{\lim}\,
        \lambda(t) g \lambda(t)^{-1} \textrm{ exists}\}$.
\item[(ii)]  $L = L_\lambda := C_G(\lambda)$.
\item[(iii)] $R_u(P) = \{g\in G \mid \underset{t\to 0}{\lim}\,
        \lambda(t) g \lambda(t)^{-1} = 1\}$.
\end{itemize}
Conversely, given  any $\lambda \in Y(G)$ 
the subset $P_\lambda$ defined as in part (i) 
is a parabolic subgroup of $G$ and
$L_\lambda$ is a Levi subgroup of $P_\lambda$. 
\end{lem}

Let $G$ act morphically on the affine algebraic variety $X$. Let $x \in X$
and let $C$ be the unique closed orbit in the closure of  $G\cdot x$.
Then there exists a subset $\Omega(x)$ of $Y(G)$
consisting of so called \emph{optimal} cocharacters $\lambda$ such 
that $\underset{t\to 0}{\lim}\,\lambda(t) \cdot x$ belongs to $C$,
\cite{kempf}, \cite{rousseau}; see \cite[\S 3]{mcninch} 
or \cite[\S 2.2]{premet} 
for the relevant parts of the theory; 
see also Slodowy's survey article \cite{slodowy2}.
We record the crucial points of this theory.

\begin{thm}
\label{thm:kempf}
Let $G$ act morphically on the affine algebraic variety $X$. Let $x \in X$
and let $\Omega(x) \subseteq Y(G)$ be the optimal class of cocharacters
defined by $x$. Then the following hold:
\begin{itemize}
\item[(i)]
$\Omega(x) \ne \varnothing$ and 
there exists an \emph{optimal} parabolic subgroup $P = P(x)$ of
$G$ so that $P = P_\lambda$ for every   
$\lambda \in \Omega(x)$.
\item[(ii)]
$\Omega(x)$ is a single $P$-orbit. 
\item[(iii)]
For every $g \in G$, we have $\Omega(g\cdot x) = g\cdot \Omega(x)$
and $P(g \cdot x) = gP(x)g\inverse$. In particular, 
$C_G(x) \le P$.
\end{itemize}
\end{thm}

\subsection{Cocharacters associated to nilpotent elements}
\label{sub:cochars}

Recall that any cocharacter $\lambda \in Y(G)$ of $G$ affords a 
$\ZZ$-grading 
\[\gg = \bigoplus_{j \in \ZZ}\gg(j, \lambda)\]
of $\gg$, where 
\[\gg(j, \lambda) := \{e \in \gg \mid \Ad(\lambda(t))e = t^je 
\text{ for every } t \in k^*\},\] cf.\ 
\cite[\S 5.5]{carter} or \cite[\S 5.1]{Jantzen}.
We recall the relevant concepts of distinguished nilpotent elements and
of cocharacters associated to a nilpotent element following 
\cite[\S 4.1, \S 5.3]{Jantzen}.

\begin{defn}
\label{def:dist}
Let $H \le G$ be a closed connected reductive 
subgroup of $G$ and $ e \in \NN \cap \hh$. 
We call $e$ \emph{distinguished} in $\hh$ if each torus
contained in $C_H(e)$ is contained in the centre of $H$.
\end{defn}

\begin{defn}
\label{def:assoc}
Let $e \in \NN$.
A cocharacter $\lambda : k^* \to G$ of $G$ is called
\emph{associated} to $e$ provided
$e \in \gg(2,\lambda)$ and 
there exists a Levi subgroup $L$ of $G$
such that $e$ is distinguished nilpotent in $\Lie L$ and 
$\lambda(k^*) \leq \DD L$.
\end{defn} 

\begin{rems}
\label{rems:1}
Let $e \in \NN$ and let 
$\lambda \in Y(G)$ that is associated to $e$. 

(i).
For $g \in C_G(e)$ the conjugate cocharacter
$g \cdot \lambda$ is also associated to $e$, cf.\ \cite[\S 5.3]{Jantzen}. 
Proposition \ref{prop:cochar1}(ii) gives a converse to this property. 

(ii).
Let $S$ be a maximal torus of $C_G(e)$. Then $e$ is distinguished 
in $\cc_\gg(S) = \Lie C_G(S)$, for $S$ is the unique maximal torus of 
$C_{C_G(S)}(e)$. Proposition \ref{prop:cochar1}(iii) 
gives a converse to this.
\end{rems}

We require some basic facts about cocharacters associated to nilpotent
elements; the following results are \cite[Rem.\ 4.7; Lem.\ 5.3]{Jantzen},
see also \cite[Prop.\ 2.5]{premet}.

\begin{prop}
\label{prop:cochar1} 
Let $e \in \NN$.
\begin{itemize} 
\item[(i)] Suppose $\Char k$ is good for $G$. Then 
cocharacters of $G$ associated to $e$ exist.
\item[(ii)] Any two cocharacters of $G$ associated to $e$ are
conjugate by an element of $C_G(e)^{\circ}$.
\item[(iii)] 
If $L$ is a Levi subgroup of $G$ with $e$ distinguished in $\Lie L$,
then the connected centre of $L$ is a maximal torus of $C_G(e)^\circ$.
\end{itemize} 
\end{prop}

\begin{rem}
\label{rem:2}
The Dynkin--Kostant classification theory giving a bijection between 
nilpotent $G$-classes and $G$-conjugacy classes of $\mathfrak{sl}(2)$-triples 
is also valid in large positive characteristic
(more precisely, when $\Char k  > 3(h-1)$, 
where $h$ denotes the Coxeter number of $G$, 
cf.\ \cite[\S\S 5.3 - 5.6]{carter}).
For $e \in \NN$, the cocharacters of $G$
constructed from the semisimple elements of 
$\mathfrak{sl}(2)$-triples containing $e$ 
(cf.\ \cite[\S 5.5]{carter})
are all associated to 
$e$ in the sense of Definition \ref{def:assoc}, see \cite[Rem.\ 5.5]{Jantzen}.
\end{rem}

Let $e \in \NN$.
In \cite[\S 2.4]{premet}, A.~Premet explicitly defines a 
cocharacter of $G$ which is associated to $e$, thanks to 
\cite[Prop.\  2.5]{premet}. 
Moreover, in  \cite[Thm.\  2.3]{premet}, Premet shows that each of these
associated cocharacters belongs to the optimal class determined by $e$.
Premet proves this under the so called \emph{standard hypotheses} on $G$, 
see \cite[\S 2.9]{Jantzen}. These restrictions were subsequently removed by 
G.\ McNinch in \cite[Prop.\ 16]{mcninch} so that this fact holds
for any reductive $G$ in good characteristic.
It thus follows from \cite[Prop.\ 16]{mcninch},
Theorem \ref{thm:kempf}, and Proposition \ref{prop:cochar1}(ii) that 
all the cocharacters of $G$ associated to $e \in \NN$ belong to the 
optimal class $\Omega(e)$ defined by $e$;
see also \cite[Prop.\ 18, Thm.\ 21]{mcninch}. 
This motivates and justifies 
the following notation which we frequently use in the sequel.
\begin{defn}
\label{d:Gamma}
Let $e \in \NN$. Then we define
\[
\Omega^a(e)  := \{\lambda \in Y(G)\mid 
\lambda  \text{ is associated to } e  \} \subseteq \Omega(e).
\]
Further, we sometimes write $\Omega_G^a(e)$ for $\Omega^a(e)$, 
to indicate that this is a set of cocharacters of $G$,
and if $H$ is a reductive subgroup of $G$ with $e \in \hh$ nilpotent 
we also write $\Omega_H^a(e)$ to denote the cocharacters of $H$ that are
associated to $e$.
\end{defn}

For $H$ a connected reductive subgroup of $G$ and $e \in \hh \cap \NN$, 
in the notation of Definition \ref{d:Gamma},
property ($\dagger$) from page \pageref{dagger} becomes the equality
\[
\Omega_H^a(e) = \Omega_G^a(e) \cap Y(H).
\]

Let $e \in \NN$ and let $\lambda \in \Omega^a(e)$. 
Let $P = P(e)$ be the canonical parabolic subgroup defined by $e$.
Then $P = C_G(\lambda) R_u(P)$ is a Levi decomposition of $P$.
Thanks to Theorem \ref{thm:kempf}(iii), $C_G(e) \le P$.
Following \cite[\S 5.10]{Jantzen} and \cite[\S 2.4]{premet}, 
we define the subgroups
\begin{equation}
\label{e:r}
C_G(e,\lambda) := C_G(e) \cap C_G(\lambda) \quad \text{ and } \quad 
R_e := C_G(e) \cap R_u(P)
\end{equation}
of $C_G(e)$.

In view of \cite[Prop.\ 2.5]{premet}, our next result is  
\cite[Thm.\ 2.3(iii)]{premet}, 
see also \cite[Prop.\ 5.10, Prop.\ 5.11]{Jantzen}.

\begin{prop}
\label{p:levidecomp}
Suppose that $\Char k$ is good for $G$.
Let $e \in \NN$ and let $\lambda \in \Omega^a(e)$.
Then $C_G(e)$ is the semidirect product of $C_G(e, \lambda)$ and $R_e$, where 
$C_G(e, \lambda)^\circ$ is reductive and $R_e = R_u(C_G(e))$.
\end{prop}

Proposition \ref{p:levidecomp} says that 
$C_G(e) = C_G(e, \lambda)R_e$ is a Levi decomposition of $C_G(e)$,
so that different choices of associated cocharacters in $\Omega^a(e)$
give conjugate Levi subgroups $C_G(e, \lambda)$ of $C_G(e)$, 
by Proposition \ref{prop:cochar1}(ii).
Our next result now readily follows from these 
two propositions.

\begin{cor}
\label{c:radical}
Suppose that $\Char k$ is good for $G$.
Let $e \in \NN$. Then $R_e$ acts simply transitively on $\Omega^a(e)$.
\end{cor}

\begin{rem}
\label{r:levi}
It follows from Propositions \ref{prop:cochar1}, \ref{p:levidecomp} 
and Corollary \ref{c:radical} that the map $\lambda \mapsto C_G(e,\lambda)$ 
is a bijection between $\Omega^a(e)$ and the set of Levi subgroups
of $C_G(e)$.
\end{rem}

Let $e \in \gg$ be nilpotent and let $\lambda \in \Omega^a(e)$. 
It follows readily from Definition \ref{def:assoc} that 
$C_G(e)$ is normalized by $\lambda(k^*)$, thus we may define the subgroup
\[
N_e := \lambda(k^*)C_G(e)
\] 
of $G$, cf.\ \cite[\S 5.3(2)]{Jantzen}.
According to Proposition \ref{prop:cochar1}(ii), $N_e$
does not depend on the choice of $\lambda \in \Omega^a(e)$;
this is also apparent, since  
$N_e =\{g \in G \mid g\cdot ke = ke\}$, see \cite[\S 2.10(2)]{Jantzen}.
Clearly, $\lambda(k^*)$ also normalizes $R_e = R_u(C_G(e))$.
Thus we may define the subgroup 
\begin{equation}
\label{e:q}
Q_e := \lambda(k^*)R_e
\end{equation}
of $N_e$.
By Corollary \ref{c:radical}, equally 
$Q_e$ does not depend on the choice of $\lambda$
in $\Omega^a(e)$.

Let $H$ be a connected reductive subgroup of $G$.
Since the nilpotent cone of $\Lie \DD H$ is a closed
subvariety of the nilpotent variety of $\hh$, and both 
are irreducible of the same dimension, we have 
$\hh \cap \NN = \Lie (\DD H) \cap \NN$.

\begin{lem}
\label{lem:derived}
Let $H$ be a connected reductive subgroup of $G$.
Let $e \in \hh \cap \NN = \Lie (\DD H) \cap \NN$ be nilpotent.
Then the cocharacters of $\DD H$ associated to $e$ are precisely  
the cocharacters of $H$ associated to $e$.
\end{lem}

\begin{proof}
Assume that $\lambda$ is a cocharacter of $\DD H$ associated to $e$.
Note that  $H = Z(H)^\circ\DD H$. 
Let $L' = C_{\DD H}(S')$ be a Levi subgroup of $\DD H$ satisfying 
the conditions of Definition \ref{def:assoc}, where $S'$ is a maximal torus of 
$C_{\DD H}(e)$.
Then $S = Z(H)^\circ S'$ is a maximal torus of $C_H(e)$ containing $S'$.
Set $L = C_H(S) = Z(H)^\circ L'$.
It follows easily that $e$ is distinguished in $\Lie L$. 
Further, $\lambda(k^*) \le \DD (L') = \DD L$, and so 
$\lambda$ is associated to $e$, viewed as a cocharacter of $H$.

The reverse implication of the lemma is shown in \cite[\S 5.6]{Jantzen}.
\end{proof}

\section{Cocharacters associated to nilpotent elements of reductive subgroups}
\label{s:cochar}

We maintain the notation and assumption of the previous sections. 
In particular, $G$ is a connected reductive algebraic group 
defined over an algebraically closed field $k$, $\Char k$ is a good prime 
for $G$, and $H$ is a closed connected reductive subgroup of $G$.

\subsection{Local Conditions}
\label{sub:local}
In this subsection we study conditions on a given nilpotent 
element $e$ in $\hh$ (or nilpotent $H$-class in $\hh$)
that ensure $(\dagger)$ from page \pageref{dagger} holds for $e$
without further assumptions on $H$ itself. 
Firstly we consider nilpotent elements $e \in \hh$ that are distinguished  
in $\gg$, Lemma \ref{lem:dist},  and secondly we study the case when 
the centralizers in $H$ and $G$ of $e$ have the same rank, 
Theorem \ref{thm:localrank}.

\begin{lem}
\label{lem:dist}
Suppose that $e \in \hh$ is nilpotent and distinguished in $\gg$. Then
$e$ is distinguished in $\hh$ and
$\Omega_H^a(e) = \Omega_G^a(e) \cap Y(H)$.
\end{lem}

\begin{proof}
Let $S$ be a torus of
$C_H(e)$. Since $e$ is distinguished in $\gg$ and $C_H(e)\leq C_G(e)$, we have
$S \leq Z(G)$. So $S \leq Z(G)\cap H\leq Z(H)$ and thus $e$ is
distinguished in $\hh$.

First let $\lambda \in \Omega_H^a(e)$.
Then $e \in \hh(2,\lambda) \subseteq \gg(2,\lambda)$.
By Lemma \ref{lem:derived}, we have 
$\lambda(k^*) \leq \DD H \leq \DD G$, and 
so $\lambda \in \Omega_G^a(e)$, as $e$ is distinguished in $\gg$.

Conversely,  let $\lambda \in \Omega_G^a(e)$ with $\lambda(k^*) \le H$. 
Since $e$ is distinguished in $\hh$, and 
since $e \in \hh \cap \gg(2,\lambda) = \hh(2,\lambda)$,
it suffices to show that $\lambda(k^*)\leq \DD H$. 
Since $e$ is distinguished in $\gg$ and  
$Z(H)^{\circ}\leq C_H(e)\leq C_G(e)$, 
we have $Z(H)^{\circ}\leq Z(G)^\circ \cap H$.
Note that $\DD H\leq \DD G \cap H \leq H$. 
Since $H$ is reductive, we
have $H = Z(H)^\circ\DD H$, so that 
$\DD G \cap H = A \DD H$, where $A$ is a
subgroup of $Z(H)^\circ$. 
By assumption, $\lambda(k^*)\leq \DD G\cap H$.
Since $\lambda(k^*)$ is connected, we have 
$\lambda(k^*)\leq (\DD G \cap H)^{\circ} = A^{\circ} \DD H$. 
Because $A^{\circ}\leq Z(H)^{\circ}\leq Z(G)^\circ \cap H$ 
and $A^{\circ}\leq \DD G\cap H$, we have
$A^{\circ}\leq \DD G \cap Z(G)^\circ$, and so 
$A^{\circ}$ is trivial.
Consequently, $\lambda(k^*)\leq \DD H$, as desired. 
\end{proof}

We give an  example for Lemma \ref{lem:dist}.

\begin{exmp}
\label{ex:dist}
Let $G$ be simple of type $E_6$ and let $H$ be the fixed point subgroup 
of the non-trivial graph automorphism of $G$; so that $H$ is of type $F_4$.
Let $C'$ be the nilpotent $H$-class with  Bala--Carter label $F_4(a_2)$
and let $C$ be the nilpotent $G$-class with  Bala--Carter label $E_6(a_3)$.
According to \cite[Table A]{lawther},  we have $C' \subset C$.
Note that each of these classes is distinguished, see \cite[\S 5.9]{carter}.
Thus Lemma \ref{lem:dist} applies.
Another example is given by the regular nilpotent 
class in $H$ which belongs to the
regular $G$-class in $\NN$; they are obviously both distinguished.
Although the results in \cite[\S 5.9]{carter} and \cite{lawther} concern
unipotent classes in $G$, since $p$ is good for $G$,
they equally apply to nilpotent classes in $\gg$.

Corollary \ref{cor:graphautomorphism} below implies that in this 
case the conclusion of Lemma \ref{lem:dist} 
holds for any nilpotent class of $H$.
\end{exmp}

\begin{lem}
\label{lem81}
Let $e \in \hh$ be nilpotent. 
Suppose that $\rank C_G(e) = \rank C_H(e)$.
Then $\Omega_H^a(e) \subseteq \Omega_G^a(e) \cap Y(H)$.
\end{lem}

\begin{proof} 
Let $\lambda \in \Omega_H^a(e)$.
Since $e \in \hh(2,\lambda) \subseteq \gg(2,\lambda)$, it
suffices to find a Levi subgroup $L$ of $G$ such that $e$ is
distinguished in $\Lie L$ and $\lambda(k^*)\leq \DD L$. 
Let $M$ be a Levi subgroup of $H$ with the properties 
as in Definition \ref{def:assoc}. 
Thanks to Proposition \ref{prop:cochar1}(iii), we have 
$M = C_H(S)$, where $S$ is a maximal torus of $C_H(e)$.
Since $C_H(e) \le C_G(e)$ and by our hypothesis, 
$S$ is also a maximal torus of $C_G(e)$. 
Thus $e$ is distinguished in $\Lie C_G(S)$, cf.\ Remark \ref{rems:1}(ii). 
Finally, since  $\lambda(k^*)\leq \DD C_H(S) \leq \DD C_G(S)$, we see that 
$\lambda \in \Omega_G^a(e)$, as desired.
\end{proof}

Let $e \in \gg$ be nilpotent. Let $\lambda$ be a cocharacter of $G$
associated to $e$. Define 
\begin{equation}
\label{e:upsilon}
\Upsilon_{\lambda}(e) =
\{S \leq C_G(e) \mid S \text{ is a maximal torus of } C_G(e) 
\text{ and } \lambda(k^*)\leq \DD C_G(S)\}.
\end{equation}
Note that by Proposition \ref{prop:cochar1}(iii),
$\Upsilon_{\lambda}(e)$ is non-empty. 

For our next result 
recall the definition of $C_G(e,\lambda) = C_G(e) \cap C_G(\lambda)$
from \eqref{e:r}.

\begin{lem}
\label{lem 10.1} 
Let $e \in \gg$ be nilpotent. Let $\lambda$ be a cocharacter of $G$
associated to $e$. Then $\Upsilon_{\lambda}(e)$ consists 
precisely of the maximal tori of $C_G(e,\lambda)$.
\end{lem} 

\begin{proof} 
Let $S\in \Upsilon_{\lambda}(e)$. Then $S$
is a maximal torus of $C_G(e)$ and 
$\lambda(k^*)\leq \DD C_G(S)\leq C_G(S)$. 
In particular, $S \leq C_G(\lambda)$ 
and so $S \leq C_G(e,\lambda)$. 
Thus $S$ is a maximal torus of $C_G(e,\lambda)$.

Conversely, let $S$ be a maximal torus of $C_G(e,\lambda)$.
Then, by what we have just shown, $S$
is conjugate in  $C_G(e,\lambda)$ to some $S' \in \Upsilon_{\lambda}(e)$,
and so in particular, $S$ is a maximal torus of $C_G(e)$.
Let $g \in C_G(e,\lambda)$ so that $S = gS'g^{-1}$. Since
$\lambda(k^*)\leq \DD C_G(S')$, we have 
$g\lambda(k^*)g^{-1} \leq g \DD C_G(S')g^{-1} = \DD C_G(S)$. 
Finally, since  $g \in C_G(\lambda)$, we obtain 
$g\lambda(k^*)g^{-1}=\lambda(k^*)$, and therefore, 
$S \in \Upsilon_{\lambda}(e)$, as desired.
\end{proof}

\begin{lem}
\label{lem:2} 
Let $e \in \hh$ be nilpotent. 
Let $\lambda \in \Omega_G^a(e)$ with $\lambda(k^*) \le H$.
Suppose there exists a maximal torus of 
$C_H(e)$ which is also a maximal torus of $C_G(e,\lambda)$. 
Then $\lambda \in \Omega_H^a(e)$. 
\end{lem}

\begin{proof}
Since $e \in \hh \cap \gg(2,\lambda) = \hh(2,\lambda)$, it
suffices to find a Levi subgroup $M$ of $H$ such that $e$ is
distinguished in $\Lie M$ and $\lambda(k^*) \leq \DD M$. Let
$L$ be a Levi subgroup of $G$ such that $e$ is distinguished in
$\Lie L$ and $\lambda(k^*)\leq \DD L$. 
By Proposition \ref{prop:cochar1}(iii), we have
$L = C_G(S)$, where $S$ is a maximal torus of $C_G(e)$. 
By Remark \ref{rems:1}(ii), $\lambda$ is a cocharacter of $C_G(S)$ that 
is associated to $e$ and $e$ is distinguished in $\cc_{\gg}(S) = \Lie C_G(S)$.
By our hypothesis and Lemma \ref{lem 10.1}, we may assume without loss 
of generality that 
$S$ is a maximal torus of $C_H(e)$.
Then $S\leq H$ and so $C_H(S)$ is a Levi subgroup of $H$.
Consequently, by Lemma \ref{lem:dist} applied to $C_H(S) \le C_G(S)$,
we get that $\lambda$ is a cocharacter of $C_H(S)$ associated to $e$.
Since $S$ is a maximal torus of $C_H(e)$ as well as of $C_{C_H(S)}(e)$, 
we can apply Lemma \ref{lem81} to $C_H(S) \le H$ and so 
$\lambda \in \Omega_H^a(e)$. 
\end{proof}

Let $e \in \hh$ be nilpotent and let $\lambda$ be a   
cocharacter of $G$ associated to $e$ with $\lambda(k^*) \le H$.
Let $P$ be the parabolic subgroup of $G$ defined by $\lambda$, i.e., 
$P = P_\lambda$, cf.\ Lemma \ref{lem:cochars}(i).
Since $\lambda$ is optimal for $e$,   
we have $P = P(e)$ is the
optimal parabolic subgroup determined by $e$, cf.\ Theorem \ref{thm:kempf}(i).
Define $P_H = P\cap H$. Since $\lambda(k^*)\leq H$, we see $P_H$ is a
parabolic subgroup of $H$. Also $P_H$ has
Lie algebra $\oplus_{i\geq 0}\hh(i,\lambda)$ and 
the unipotent radical $R_u(P_H)$ of $P_H$ is
$R_u(P)\cap H$ and $C_H(\lambda) = C_G(\lambda)\cap H$ is a
Levi subgroup of $P_H$. The unipotent radical of $P_H$ and the
Levi subgroup $C_H(\lambda)$ have Lie algebras
$\oplus_{i>0}\hh(i,\lambda)$ and $\hh(0,\lambda)$,
respectively, see \cite[\S 5.1]{Jantzen}.
Analogous to \eqref{e:r} define the subgroups 
\begin{equation}
\label{eq:c_H}
C_H(e,\lambda) := C_H(e)\cap C_H(\lambda)
\quad \text{ and }\quad U_e := C_H(e) \cap R_u(P_H)
\end{equation}
of $C_H(e)$.
Note that $C_H(e,\lambda)\cap U_e=\{1\}$, so $C_H(e,\lambda)U_e$
is a semidirect product of algebraic groups. Next we require the
following facts, see \cite[\S5.10]{Jantzen}. For all $i \in \ZZ$ we
have
\begin{equation}
\label{eq1} 
\Ad(x)(\hh(i,\lambda))=\hh(i,\lambda)
\text{ for all } x \in C_H(\lambda);
\end{equation}
\begin{equation}
\label{eq2}
(\Ad(y)-1)(\hh(i,\lambda))\subseteq \bigoplus_{j> i}\hh(j,\lambda)
\text{ for all } y \in R_u(P_H).
\end{equation}
\goodbreak

\begin{lem}
\label{lem 10.2} 
With the notation introduced in \eqref{eq:c_H}, we have the following.
\begin{itemize}
\item[(i)]
$C_H(e) \le P_H$.
\item[(ii)]
$U_e$ is a normal connected unipotent subgroup of $C_H(e)$.
\item[(iii)]
$C_H(e)$ is the semidirect product of $C_H(e,\lambda)$ and $U_e$.
\end{itemize}
\end{lem} 

\begin{proof}
(i). 
Let $x \in C_H(e) = C_G(e)\cap H$. 
Since $\lambda \in \Omega^a(e)$, we have $C_G(e)\subseteq P$,
by Theorem \ref{thm:kempf}(iii).
Therefore, $x \in P\cap H = P_H$.

(ii).
By definition, $U_e$ is unipotent.  
Let $x \in U_e$ and $y \in C_H(e)$. Now $x \in
R_u(P_H)$ and $y \in P_H$. So $yxy^{-1} \in R_u(P_H)$, as
$R_u(P_H)$ is normal in $P_H$. Also $x, y \in C_H(e)$ so $yxy^{-1} \in
R_u(P_H)\cap C_H(e) = U_e$.  
Since $U_e$ is normalized by $\lambda(k^*)$, it is connected:
For any $x\in U_e$ we
have a morphism $\phi_x : k \to U_e$ given by
$\phi_x(t) = \lambda(t) x \lambda(t)\inverse$ for $t\in k^*$
and $\phi_x(0) = \underset{t\to 0}{\lim}\,\lambda(t) x \lambda(t)\inverse = 1$,
cf.\ Lemma \ref{lem:cochars}(iii).
Thus for each $x\in U_e$ 
the image of $\phi_x$ is a connected subvariety of $U_e$ 
containing $1$ and $x = \phi_x(1)$.  

(iii).
Let $z \in C_H(e)$. Thanks to part (i), we can write $z$ 
as $z = xy$ with $x \in C_H(\lambda)$ and $y \in R_u(P_H)$. 
Since $e \in \hh(2,\lambda)$, it follows from \eqref{eq2} that
$\Ad(y)(e)=e+e'$ with $e' \in \oplus_{i\geq 3}\hh(i,\lambda)$. 
Since $z \in C_H(e)$, we have $\Ad(x)(e+e') = e$. It follows 
that $\Ad(x)(e') = 0$, thanks to \eqref{eq1}. 
As $\Ad(x) : \hh \to \hh$ is a Lie algebra automorphism of $\hh$, 
we infer that $e' = 0$. 
Consequently, $y \in C_H(e) \cap R_u(P_H)$ and thus 
$x \in C_H(e)\cap C_H(\lambda)$.
Therefore, $C_H(e) = C_H(e,\lambda)U_e$.
\end{proof}

\begin{rem} 
Despite the analogy between Proposition \ref{p:levidecomp} and 
Lemma \ref{lem 10.2}, the semidirect product 
$C_H(e) = C_H(e,\lambda)U_e$ need not be a Levi decomposition of $C_H(e)$;
we cannot invoke Proposition \ref{p:levidecomp}, as we do not know whether
$\lambda$ lies in $\Omega_H^a(e)$.
It follows from 
Theorem \ref{thm:main} that this is the case 
for the  particular class of subgroups $H$ of $G$ considered there.
\end{rem}

\begin{lem}
\label{lem 10.45}
Let $e \in \hh$ be nilpotent. 
Suppose that $\rank C_G(e) = \rank C_H(e)$. 
Let $\lambda  \in \Omega_G^a(e) \cap Y(H)$.
Then there exists a maximal torus of
$C_H(e)$ which is also a maximal torus of $C_G(e,\lambda)$.
\end{lem}

\begin{proof} 
It follows from Lemma \ref{lem 10.2}(iii) that 
$\rank C_H(e) = \rank C_H(e, \lambda)$. Obviously, we have 
$C_H(e, \lambda) \leq C_G(e, \lambda) \le C_G(e)$. 
It follows from the assumption on the rank of the centralizers 
of $e$ in $H$ and in $G$ that a maximal torus of $C_H(e,\lambda)$, 
therefore of $C_H(e)$, is also a maximal torus of $C_G(e,\lambda)$.
\end{proof} 

Our next result is a consequence of Lemmas \ref{lem:2} 
and \ref{lem 10.45}.

\begin{prop}
\label{prop2}
Let $e \in \hh$ be nilpotent. Suppose that $\rank C_G(e) = \rank C_H(e)$. 
Then $\Omega_G^a(e) \cap Y(H) \subseteq \Omega_H^a(e)$.
\end{prop}

Our next result follows immediately 
from Lemma \ref{lem81} and Proposition \ref{prop2}.

\begin{thm}
\label{thm:localrank} 
Let $e \in \hh$ be nilpotent. Suppose that $\rank C_G(e) = \rank C_H(e)$. 
Then $\Omega_H^a(e) = \Omega_G^a(e) \cap Y(H)$.
\end{thm}

We give an example for Theorem \ref{thm:localrank}. 

\begin{exmp}
\label{ex:1}
We return to the case of Example \ref{ex:dist}.
Let $G$ be simple of type $E_6$ and let $H$ be the 
standard subgroup 
of $G$ of type $F_4$. 
Let $D'$ be the nilpotent $H$-class with Bala--Carter label $\tilde A_2$
and let $D$ be the nilpotent $G$-class with Bala--Carter label $2A_2$.
According to \cite[Table A]{lawther}, we have $D' \subset D$.
It is known that the reductive parts of the centralizers  
of these classes in $H$ and $G$ are of type $G_2$ in each case, 
see \cite[Ch.\ 13]{carter}.
In particular, the rank condition in Theorem \ref{thm:localrank}
is satisfied for an element belonging to $D'$.
This also applies to the classes with Bala--Carter labels $C_3$ in $H$ and
$A_5$ in $G$, here the reductive parts of the centralizers 
are of type $A_1$ in each case. 
There is one further pair of classes with the same property.

It follows from Corollary \ref{cor:graphautomorphism} below that
in this example  
the conclusion of Theorem \ref{thm:localrank} 
holds for any nilpotent element of $\hh$ irrespective of the 
condition on the ranks of the respective centralizers.
\end{exmp}

See also Example \ref{ex:g2f4} below for another application of 
Theorem \ref{thm:localrank}. 

\subsection{Global Conditions}
\label{sub:global}
We maintain the notation from the previous sections.
In this subsection we study conditions on the subgroup $H$ of $G$ 
(rather than on a given nilpotent $H$-class in $\hh$) 
that ensure that $(\dagger)$ holds for all $e \in \NN \cap \hh$.

Our first result shows that the reverse inclusion of 
$(\dagger)$ always holds in good characteristic 
without any restrictions on $H$.

\begin{prop}
\label{eqweak} 
Let $e \in \hh$ be nilpotent. 
Then $\Omega_G^a(e) \cap Y(H) \subseteq \Omega_H^a(e)$.
\end{prop}

\begin{proof} 
Let $\lambda \in \Omega_G^a(e)$ with $\lambda(k^*) \le H$.
First assume that $e$ is distinguished in $\hh$. 
Since $e \in \hh \cap \gg(2,\lambda) = \hh(2,\lambda)$ 
and $e$ is distinguished in $\hh$, it
suffices to show that $\lambda(k^*)\leq \DD H$. 
Set $Z = Z(H)^{\circ}$. Then $Z \leq C_H(e,\lambda) \leq C_G(e,\lambda)$. 
Choose a maximal torus $S$ of $C_G(e,\lambda)$ so that $Z \leq S$. 
Thanks to Lemma \ref{lem 10.1}, $S$ is a maximal torus of $C_G(e)$.
Since $S \leq C_G(e,\lambda)$, we have $\lambda(k^*) \leq \DD C_G(S)$, 
by Lemma \ref{lem 10.1}.
Also, as $Z \leq S$, we have 
$C_G(S) \leq C_G(Z)$. Thus $\lambda(k^*) \leq \DD C_G(Z)$ and so
$\lambda(k^*) \leq H \cap \DD C_G(Z)$. 
As $H \leq C_G(Z)$, we get $\DD H \leq H \cap \DD C_G(Z) \leq H$. Since
$H$ is reductive, we have $H = Z \DD H$, 
so $H \cap \DD C_G(Z) = A\DD H$, 
where $A \leq Z$. As $\lambda(k^*)$ is
connected, $\lambda(k^*) \leq A^\circ \DD H$. 
Clearly, we have $A^{\circ} \leq A \DD H = H \cap \DD C_G(Z)$ 
and $A^{\circ} \leq Z$ so $A^{\circ}\leq Z \cap \DD C_G(Z)$. 
Since $C_G(Z)$ is a connected reductive subgroup of $G$ 
and $Z$ is contained in the connected centre of $C_G(Z)$, it follows that  
$Z \cap \DD C_G(Z)$ is finite. Thus
$A^{\circ}$ is trivial and hence $\lambda(k^*) \leq \DD H$, as desired.

Now we consider the general case where $e$ is not necessarily  
distinguished in $\hh$.
Let $S$ be a maximal torus of $C_H(e,\lambda)$.
By Lemma \ref{lem 10.2}(iii), $S$ is then also a maximal torus of $C_H(e)$.
Since $S \le C_G(\lambda)$, we have $\lambda(k^*)\leq C_H(S)$ 
and that $e$ is distinguished in $\Lie C_H(S)$, 
cf.\ Remark \ref{rems:1}(ii). Thus, by
the distinguished case just proved, 
$\lambda$ is a
cocharacter of $C_H(S)$ associated to $e \in \cc_\hh(S) = \Lie C_H(S)$. 
Observe that $\dim S = \rank C_{C_H(S)}(e) = \rank C_H(e)$. 
It thus follows from Lemma \ref{lem81} applied
to $C_H(S) \le H$ that $\lambda$ is a 
cocharacter of $H$ associated to $e$, as claimed.
\end{proof}

Our next result shows that the forward inclusion of 
$(\dagger)$ holds provided there is at least one cocharacter of $G$
that is associated to $e \in \NN \cap \hh$ and takes values in $H$.

\begin{lem}
\label{lem 11.5} 
Let $e \in \hh \cap \NN$.
Then $\Omega_H^a(e)\subseteq \Omega_G^a(e) \cap Y(H)$,
provided $\Omega_G^a(e) \cap Y(H)$ is non-empty. 
\end{lem}

\begin{proof} 
Let $\lambda \in \Omega_G^a(e) \cap Y(H)$ and let $\mu \in \Omega_H^a(e)$. 
Then,  by Proposition \ref{eqweak}, $\lambda$ is a cocharacter of $H$ 
associated to $e$. 
So, by Proposition \ref{prop:cochar1}(ii), there exists an 
$x \in C_H(e)^{\circ}$ such that $x\lambda x^{-1} = \mu$. Clearly,
$C_H(e)^{\circ} \leq C_G(e)^{\circ}$. By Remark \ref{rems:1}(i)
it follows that $\mu$ is thus a cocharacter of $G$ associated to $e$.
\end{proof}

Combining Proposition \ref{eqweak} and Lemma \ref{lem 11.5}
gives our next result.

\begin{thm}
\label{thm:reduction}
Let $e\in \hh \cap \NN$.
If $\Omega_G^a(e) \cap Y(H)$ is non-empty, 
then $\Omega_H^a(e) = \Omega_G^a(e) \cap Y(H)$.
\end{thm}

Let $S$ be a linearly reductive group acting on $G$ by automorphisms, so that  
$G$ is an $S$-group.
Set $H = C_G(S)^\circ$ which is reductive, by 
Proposition \ref{prop:linear1}(i).
Let $e \in \hh = \cc_\gg(S)$ be nilpotent. 
In the proof of \cite[Thm.\ C]{rich0}, R.W.~Richardson showed that there 
always exists a cocharacter of $G$ which belongs to
the optimal class $\Omega(e)$ defined by $e$ 
so that its image lies in $H$, see \cite[\S 8]{rich0}.
For our purpose we need a variant of this result:
we require the existence of 
a cocharacter of $G$ which is associated to $e$
(rather than merely being optimal)
so that its image lies in $H$. 
That is, we need to construct an $S$-fixed cocharacter 
of $G$ which is associated to $e$. This is done in 
Lemma \ref{l:S-fixedtau} with the aid of Richardson's result
Proposition \ref{prop:richardson}.
For our next result, recall the subgroup $Q_e$ of $N_e$ from 
\eqref{e:q}. 

\begin{lem}
\label{l:S-stablecentralizer}
Let $e \in \hh = \cc_\gg(S)$ be nilpotent. Then 
\begin{itemize}
\item[(i)]
$\Omega^a(e)$ is $S$-stable;
\item[(ii)]
$C_G(e)$ is $S$-stable;
\item[(iii)]
$Q_e$ is $S$-stable.
\end{itemize}
\end{lem}

\begin{proof}
(i).
Let $\lambda \in \Omega^a(e)$.
Since the induced actions of $S$ and $G$ on $\gg$ are
compatible in the sense of \eqref{e:compatible} and linear,
we have
\[
\Ad((s\cdot \lambda)(t))e = s \cdot \Ad(\lambda(t))e = s\cdot (t^2e) = t^2 e,
\]
for every $s \in S$ and $t \in k^*$; so $e \in \gg(2, s \cdot \lambda)$
for every $s \in S$.
Clearly, if $e$ is distinguished in $\Lie L$, then $e = s\cdot e$ is 
distinguished in $\Lie (s\cdot L)$ and $s\cdot L$ is another Levi subgroup of
$G$ for $s \in S$.
Finally, $\lambda(k^*) \le \DD L$ implies that 
$(s\cdot \lambda)(k^*) \le s\cdot \DD L = \DD(s\cdot L)$ for $s \in S$.
It follows that $s \cdot \lambda \in \Omega^a(e)$ for any $s \in S$.

(ii).
Again, by the compatibility of the actions of $S$ and $G$ on $\gg$, 
we obtain
\[
\Ad(s \cdot g)e =  s\cdot\Ad(g) e = s \cdot e = e,
\]
for any $g \in C_G(e)$, $s \in S$.

(iii).
Since $C_G(e)$ is $S$-stable, so is $R_e = R_u(C_G(e))$.
Since $\Omega^a(e)$ is $S$-stable, by part (i), 
it follows from Corollary \ref{c:radical} that
$Q_e$ is also $S$-stable.
\end{proof}

\begin{lem}
\label{l:S-fixedtau}
Let $e \in \hh = \cc_\gg(S)$ be nilpotent. Then the following hold.
\begin{itemize}
\item[(i)] 
There exists an $S$-stable Levi subgroup $\lambda(k^*)$ of $Q_e$.
\item[(ii)] 
There exists an $S$-fixed cocharacter in $\Omega^a(e)$.
\item[(iii)] 
There exists an $S$-stable Levi subgroup
$C_G(e, \lambda)$ of $C_G(e)$.
\end{itemize}
\end{lem}

\begin{proof}
(i).
Let $\lambda' \in \Omega^a(e)$ so that $Q_e = \lambda'(k^*)R_e$.
Then $R_u(Q_e) = R_e$.
According to Corollary \ref{c:radical}, $R_e$ acts simply transitively 
on $\Omega^a(e)$. Thus $R_e$ acts simply transitively on the set of
Levi subgroups of $Q_e$, Remark \ref{r:levi}. 
Thanks to Lemma \ref{l:S-stablecentralizer}(iii),  $Q_e$ is $S$-stable.
The desired result now follows from Proposition \ref{prop:richardson}.

(ii).
Let $\lambda(k^*)$ be as in (i), i.e., 
$(s \cdot \lambda)(k^*) = s \cdot (\lambda(k^*)) 
= \lambda(k^*)$ for every $s \in S$.
Since $e \in \gg(2, \lambda) \cap \gg(2, s \cdot \lambda)$,
it follows from \cite[Lem.\ 4.11]{Jantzen} that 
$s \cdot \lambda = \lambda$ for every $s \in S$. 

(iii).
Let $\lambda$ in $\Omega^a(e)$ be $S$-fixed as in (ii).
Since $S$ and $G$ act compatibly on $Y(G)$ 
in the sense of \eqref{e:compatible}, $C_G(\lambda)$ is $S$-stable. For,
\[
(s\cdot g)\cdot \lambda(t)= s \cdot (g \cdot \lambda)(t) 
= s \cdot \lambda(t) = \lambda(t),
\]
for all $s \in S$, $g \in C_G(\lambda)$, and $t \in k^*$.
Consequently, since $C_G(e)$ is $S$-stable, 
by Lemma \ref{l:S-stablecentralizer}(ii),
so is $C(e,\lambda) = C_G(e) \cap C_G(\lambda)$.
The result now follows from Proposition \ref{p:levidecomp}.
\end{proof}

Finally, Theorem \ref{thm:main0} follows from 
Theorem \ref{thm:reduction} and Lemma \ref{l:S-fixedtau}(ii):

\begin{thm}
\label{thm:main}
Let $S$ be a linearly
reductive group acting on $G$ by automorphisms and 
set $H = C_G(S)^\circ$. 
Let $e\in \hh \cap \NN$.
Then $\Omega_H^a(e) = \Omega_G^a(e) \cap Y(H)$.
\end{thm}

We record various special cases of Theorem \ref{thm:main} 
as separate corollaries.

Since a Levi subgroup of $G$ is of the form $C_G(S)$ for some 
torus $S$ of $G$, our next result is immediate from 
Theorem \ref{thm:main}.

\begin{cor}
\label{cor:levi}
Let $H$ be a Levi subgroup of $G$.
Let $e\in \hh \cap \NN$.
Then $\Omega_H^a(e) = \Omega_G^a(e) \cap Y(H)$.
\end{cor}

Likewise, our next result is immediate from Theorem \ref{thm:main}
and Remark \ref{rem:ssauto}.

\begin{cor}
\label{cor:sautomorphism}
Let $\sigma$ be a semisimple automorphism of $G$.
Let $H = C_G(\sigma)^\circ$.  
Let $e\in \hh \cap \NN$.
Then $\Omega_H^a(e) = \Omega_G^a(e) \cap Y(H)$.
\end{cor}

Here is a special case of Corollary \ref{cor:sautomorphism}.
Assume that $G$ is simple and let 
$\gamma$ be a non-trivial graph automorphism of $G$. 
If $\Char k$ is coprime to the order of $\gamma$, 
then $\gamma$ is a semisimple automorphism of $G$.

\begin{cor}
\label{cor:graphautomorphism}
Let $G$ be simple and let 
$\gamma$ be a non-trivial graph automorphism of $G$. Suppose that
$\Char k$ is coprime to the order of $\gamma$.
Let $H = C_G(\gamma)^\circ$. 
Let $e\in \hh \cap \NN$.
Then $\Omega_H^a(e) = \Omega_G^a(e) \cap Y(H)$.
\end{cor}

We give two examples for Corollary \ref{cor:graphautomorphism}. 

\begin{exmp}
\label{ex:graph1}
Suppose that $p \ne 2$.
Let $V$ be a finite-dimensional $k$-vector space.
Let $H$ be either $\SP(V)$ or $\SO(V)$.
Observe that $H$ is the fixed point subgroup of
an involution of $\SL(V)$ (cf.\ \cite[\S 11 p.\  169]{St1})
and thus Corollary \ref{cor:graphautomorphism}
applies. That is, 
the cocharacters of the classical groups $H = \SP(V)$ or $\SO(V)$ 
associated to a given nilpotent element $e$ in the Lie algebra of $H$
are precisely the cocharacters of the ambient linear group
$\SL(V)$ associated to $e$ whose image lies in $H$.
\end{exmp}

\begin{exmp}
\label{ex:graph2}
Suppose that $\Char k > 3$. Let $G$ be of type $D_4$ and let $\gamma$ be
the triality graph automorphism of $G$. Then $H = C_G(\gamma)^\circ$ is 
of type $G_2$ and so Corollary \ref{cor:graphautomorphism} applies.
Thus, for a given nilpotent element $e \in \hh$ 
we can realize the cocharacters of $H$ associated to $e$ 
as the cocharacters of $G$ associated to $e$ that are $\gamma$-fixed.
\end{exmp}

Recall from subsection \ref{sub:regular} that  
for $s \in G$ semisimple, $C_G(s)^{\circ}$ 
is a pseudo-Levi subgroup of $G$. The next result 
is again a special case of Corollary \ref{cor:sautomorphism}.

\begin{cor}
\label{cor:pseudolevi}
Let $H$ be a pseudo-Levi subgroup of $G$.
Let $e\in \hh \cap \NN$.
Then $\Omega_H^a(e) = \Omega_G^a(e) \cap Y(H)$.
\end{cor}

Clearly, each pseudo-Levi subgroup of $G$ is of maximal rank.
The subsystems corresponding to maximal rank, semisimple 
subgroups  of a simple group $G$ are determined by means of the algorithm of 
Borel and de Siebenthal \cite{BoSe}, see also 
\cite[Ex.\ Ch.\ VI \S 4.4]{bou}.
Using Corollary \ref{cor:pseudolevi}, the algorithm of 
Borel and de Siebenthal, as well as Deriziotis'
characterization of maximal rank reductive subgroups 
(cf.\ \cite[\S 2.15]{Humphreys}), 
we can generalize Corollary \ref{cor:pseudolevi} to arbitrary
maximal rank reductive subgroups.

Our next result is Theorem \ref{thm:maxrank0}.

\begin{thm}
\label{thm:maxrank}
Let $H$ be a connected reductive subgroup of $G$ of maximal rank.
Let $e\in \hh \cap \NN$.
Then $\Omega_H^a(e) = \Omega_G^a(e) \cap Y(H)$.
\end{thm}

\begin{proof}
Thanks to Lemma \ref{lem:derived} and 
by passing to simple factors, we may assume that $G$ is simple. 
Let $H$ be a maximal rank reductive subgroup of $G$. Let $T$ be 
a maximal torus of $G$ contained in $H$.
Let $\Pi$ be a set of simple roots of $\Psi = \Psi(G,T)$, let $\varrho$
be the highest root of $\Psi^+$, and let $W$ be the Weyl group of $G$ 
with respect to $T$. Let $\Phi = \Phi(H,T)$ be the root system of $H$; 
in particular $\Phi$ is a semisimple subsystem of $\Psi$. 
Since  $\Char k$ is good for $G$, it follows from Deriziotis' Criterion
(cf.\ \cite[\S 2.15]{Humphreys}) that 
$H$ is of the form $H = C_G(s)^\circ$ 
if and only if $\Phi$ admits a base which is $W$-conjugate
to a proper subset of $\Pi \cup \{-\varrho\}$.
This construction coincides with the inductive step in the
Borel--de Siebenthal procedure, \cite{BoSe}.
Thanks to \cite[Prop.\ 16]{mcninchsommers}, 
since $p$ is good for $G$, it is also good for $H$.
Because every maximal rank subsystem of $\Psi$ 
is obtained by an iteration of the 
Borel--de Siebenthal procedure, the proposition follows by
a repeated application of 
this algorithm, Deriziotis' Criterion,
and Corollary \ref{cor:pseudolevi}.
\end{proof}

\begin{thm}
\label{thm:regular}
Let $H$ be a connected regular reductive subgroup of $G$.
Let $e\in \hh \cap \NN$.
Then $\Omega_H^a(e) = \Omega_G^a(e) \cap Y(H)$.
\end{thm}

\begin{proof}
Let $T$ be a maximal torus of $G$ normalizing $H$. Then $TH$ is a connected 
reductive subgroup of $G$ of maximal rank.  
Since $\DD(TH) = \DD H$ is the semisimple part of $H$
and so $TH = Z(TH)\DD H$, the result follows from 
Lemma \ref{lem:derived} and Theorem \ref{thm:maxrank}.
\end{proof}

Our next two examples show that 
by iterating our results we can ensure that $(\dagger)$ holds
even in cases where we cannot apply Theorem \ref{thm:main} directly.

\begin{exmp}
\label{ex:2}
Let $G$ be of type $F_4$ and let $K$ be a connected simple subgroup of $G$ of
type $D_4$ and $H$ a connected subgroup of $K$ of type $G_2$.
Suppose that $p$ is good for $G$.
Then a successive application of Theorem \ref{thm:maxrank} and 
the conclusion from Example \ref{ex:graph2} show that 
$(\dagger)$ also holds for the embedding $H \le G$.
Note that $H$ is not a regular subgroup of $G$ so that we cannot 
invoke Theorem \ref{thm:regular} directly to  the embedding $H \le G$.
This is the \emph{standard} embedding of $G_2$ in $F_4$.
We discuss a different embedding of $G_2$ in $F_4$ 
in characteristic $7$ in Example \ref{ex:g2f4} below.
\end{exmp}

\begin{exmp}
\label{ex:3}
Let $G$ be of type $E_8$ and suppose that $p$ is good for $G$. Let
$K$ be a maximal rank subgroup of $G$ of type $D_4\times D_4$, and let
$H$ be a connected subgroup of $K$ of type $D_4$ embedded diagonally into $K$.
Then $H$ is the fixed point subgroup of the involution interchanging the 
$D_4$-factors in $K$. Note that $H$ is not a regular subgroup of $G$,
so we cannot invoke  Theorem \ref{thm:regular}.
Nevertheless, it follows from Corollary \ref{cor:sautomorphism} and 
Theorem \ref{thm:maxrank} that $(\dagger)$ 
holds for the embedding $H \le G$.
\end{exmp}

\begin{rems}
(i). Theorem \ref{thm:maxrank} answers 
a question  posed by J.C.~Jantzen, cf.\ \cite[\S 5.12]{Jantzen}.

(ii).
The forward inclusion of Corollary \ref{cor:pseudolevi} was already proved  
in \cite[Prop.\ 23, Rem.\ 25]{mcninchsommers} by different methods. 
\end{rems}

\begin{exmp}
\label{ex:g2f4}
In \cite[Thm.\ 1(c)]{testerman}, D.\ Testerman showed that 
there is a maximal subgroup of type $G_2$ in $F_4$ in 
characteristic $p = 7$; see also \cite[Thm.\ 1]{seitz}.
In this case let $G$ be the ambient group 
of type $F_4$ and let $H$ be the maximal subgroup of type $G_2$.
The fusion of the unipotent classes 
of this embedding has been determined 
by R.\ Lawther (unpublished). 
Since $p = 7$ is good for $G$,
this also determines the fusion of the nilpotent classes of $\hh$ in $\gg$;
this is given in terms of the corresponding Bala--Carter labellings
as follows:
\goodbreak

\begin{table}[h]
\renewcommand{\arraystretch}{1.5}
\begin{tabular}{|c|ccccc|}
\hline
$G_2$ & 1 & $A_1$ & $\tilde A_1$ & $G_2(a_1)$ & $G_2$ \\
\hline
$F_4$ & 1 & $A_1\tilde A_1$ & $\tilde A_2 A_1$ & $F_4(a_3)$ & $F_4(a_1)$\\
\hline
\end{tabular}
\bigskip
\caption{The fusion of nilpotent classes for $G_2 \le_{\max} F_4$ ($p=7$).} 
\label{t:1}
\end{table}


In contrast to the standard (non-maximal) embedding of $G_2$ in $F_4$ which
exists in any characteristic (Example \ref{ex:2}), we cannot deduce 
$(\dagger)$ for the embedding $H \le G$ directly or by iterating our results, 
as $H$ is maximal in $G$.
Nevertheless, our methods allow us to easily deduce 
that $(\dagger)$ holds for all but one of the nilpotent classes of $\hh$. 

Note that the trivial cocharacter 
is associated to $e=0$ for both $H$ and $G$.
If $e$ belongs to the regular or subregular class in $\hh$, then according to 
Table \ref{t:1} and \cite[\S 5.9]{carter}, the corresponding $G$-classes 
$F_4(a_3)$ and  $F_4(a_1)$ in $\NN$
are distinguished. Therefore, the result follows by Lemma \ref{lem:dist}.
Let $e$ belong to the $G_2$-class with label $\tilde A_1$.
According to Table \ref{t:1} and 
\cite[Ch.\ 13]{carter}, the reductive parts of the 
centralizers $C_H(e)$ and $C_G(e)$
of this class in $G_2$ and the corresponding class 
$\tilde A_2 A_1$ in  $F_4$ are both of type $A_1$. 
In particular, $C_H(e)$ and $C_G(e)$
have the same rank and Theorem \ref{thm:localrank} gives the desired 
result in this case.

For the remaining pair of nilpotent orbits one can
use Theorem \ref{thm:reduction} and show directly that $(\dagger)$ also holds
in this case. 
Fix maximal tori $T_H$ and $T$ of $H$ and $G$ respectively, so that
$T_H \le T$. Let $\alpha$ be the long simple root of $H$ with respect to $T_H$
and let $\{\alpha_1, \ldots, \alpha_4\}$ be the set of simple roots of $G$
with respect to $T$ so that $\alpha_1$ and $\alpha_2$ are long.
Let $e = e_\alpha$, a non-trivial root vector 
in the root space of $\alpha$ in $\gg$. 
Then $e$ belongs to the $H$-class in $\hh$ with 
label $A_1$ and the coroot $\alpha \check \in Y(H)$ 
is an associated cocharacter
of $e$, see \cite[\S 5.13]{Jantzen}.
Let $L_{\{\alpha_1, \alpha_3\}}$ be the standard Levi subgroup 
of $G$ of type $A_1\tilde A_1$ with root system $\{\pm\alpha_1, \pm\alpha_3\}$.
Let $M = {}^{s_{\alpha_2}}(L_{\{\alpha_1, \alpha_3\}})$, where
$s_{\alpha_2}$ is the simple reflection in the Weyl group of $G$ 
corresponding to $\alpha_2$. Then $M$ is another Levi subgroup of $G$. 
Using Lawther's explicit fusion calculations, one can show that 
$\alpha \check(k^*) \le \DD M$ and that 
$e$ is distinguished in $\mm$.
Finally, since $e \in \hh(2,\alpha \check) \subset \gg(2,\alpha \check)$, 
it follows that $\alpha \check \in \Omega_G^a(e)$.
Therefore, $\alpha \check \in \Omega_G^a(e)\cap Y(H)$, and consequently, 
by Theorem \ref{thm:reduction}, we have 
$\Omega_H^a(e) = \Omega_G^a(e) \cap Y(H)$.
\end{exmp}


\bigskip 
{\bf Acknowledgements}:
The first author acknowledges funding by the EPSRC. 
We are grateful to R.\ Lawther for making available 
information on the fusion of the nilpotent classes given in Table \ref{t:1}
in Example \ref{ex:g2f4} from an unpublished manuscript. 
Further, we would like to thank S.\ Goodwin for carefully reading 
earlier versions of the paper and for suggesting various improvements.
We are particularly indebted to G.\ McNinch for 
very helpful discussions and comments.

\bigskip



\begin{thebibliography}{00}

\bibitem{BaRi} 
P.~Bardsley, R.W.~Richardson, 
\emph{\'Etale slices for algebraic transformation 
groups in characteristic $p$}, 
Proc.\  London Math.\  Soc.\  (3) \textbf{51} (1985), no. 2, 295--317. 

\bibitem{borel} 
A.~Borel, \emph{Linear Algebraic Groups}, Graduate Texts in 
Mathematics, \textbf{126}, Springer-Verlag 1991.

\bibitem{BoSe}
A.~Borel, J.~de Siebenthal, 
\emph{Les sous-groupes ferm\'es de rang maximum des groupes de Lie clos},
Comment.\  Math.\  Helvet.\  \textbf{23}, (1949), 200--221.

\bibitem{bou} 
N.~Bourbaki, 
\emph{Groupes et alg\`{e}bres de Lie}, 
Chapitres 4, 5 et 6, Hermann, Paris, 1975. 

\bibitem{carter}
R.~W.~Carter,
\emph{Finite groups of Lie type. Conjugacy classes and
complex characters}, Pure and Applied Mathematics,
New York, 1985.

\bibitem{Humphreys} 
J.E.~Humphreys, 
\emph{Conjugacy classes in semisimple algebraic groups}; 
Mathematical Surveys and Monographs, 43. 
American Mathematical Society, Providence, RI, 1995. 

\bibitem{Jantzen} 
J.C.~Jantzen,
\emph{Nilpotent Orbits in Representation Theory}, 
in Lie Theory. Lie Algebras and Representations.
Progress in Math. vol. 228,
J.-P. Anker, B. Orsted, eds.
Birkh\"auser Boston, 2004.

\bibitem{kempf}
G.R.~Kempf, 
\emph{Instability in Invariant Theory}, 
Ann.\  Math.\ \textbf{108} (1978), 299--316.

\bibitem{lawther}
R.~Lawther, 
\emph{Jordan block sizes of unipotent elements in exceptional 
algebraic groups}, Comm. Algebra \textbf{23} (1995), no. 11, 4125--4156.

\bibitem{mcninch}
G.~McNinch, 
\emph{Nilpotent orbits over ground fields of good characteristic}, 
Math. Ann. \textbf{329} (2004), no. 1, 49--85. 

\bibitem{mcninch2}
\bysame, 
\emph{Optimal ${\rm SL}(2)$-homomorphisms}, 
Comment. Math. Helv. \textbf{80} (2005), no. 2, 391--426. 

\bibitem{mcninchsommers}
G.~McNinch, E.~Sommers,
\emph{Component groups of unipotent centralizers in good characteristic}, 
Special issue celebrating the 80th birthday of Robert Steinberg. 
J. Algebra \textbf{260} (2003), no. 1, 323--337.

\bibitem{mostow}
G.D.~Mostow,
\emph{Fully reducible subgroups of algebraic groups},
Amer. Math J., \textbf{78}, (1956), 200--221.

\bibitem{nagata}
M.~Nagata,
\emph{Complete reducibility of rational representations of a matric group},
J.\ Math.\ Kyoto University \textbf{1} (1961), 87--99.

\bibitem{premet2} 
A.~Premet, 
\emph{An analogue of the Jacobson-Morozov Theorem for Lie algebras of reductive
groups of good characteristics}, 
Trans. Amer. Math. Soc. \textbf{347}, (1995), 2961--2988.

\bibitem{premet} 
\bysame, 
\emph{Nilpotent orbits in good characteristic 
and the Kempf-Rousseau theory},
Special issue celebrating the 80th birthday of Robert Steinberg. 
J. Algebra \textbf{260} (2003), no. 1, 338--366.

\bibitem{rich0}
R.W. ~Richardson,  
\emph{On orbits of algebraic groups and Lie groups},
Bull.\  Austral.\  Math.\  Soc.\  \textbf{25} (1982), no.\  1, 1--28.

\bibitem{rousseau}
G.~Rousseau,
\emph{Immeubles sph\'eriques et th\'eorie des invariants},
C. R. Acad. Sci. Paris \textbf{286} (1987), 247--250.

\bibitem{seitz}
G.M.~Seitz, 
\emph{Maximal subgroups of exceptional algebraic groups}, 
Mem. Amer. Math. Soc. \textbf{441}, (1991).

\bibitem{seitz2}
\bysame, 
\emph{Unipotent elements, tilting modules, and saturation},
Invent. Math. \textbf{141} (2000), no. 3, 467--502.

\bibitem{slodowy2}
P.~Slodowy, 
\emph{Theorie der optimalen Einparameteruntergruppen},
Algebraische Trans\-formations\-gruppen und Invariantentheorie
(eds.\  H.\  Kraft, P.\  Slodowy and T.A.\  Springer),
DMV Seminar, \textbf{13}, Birkh\"auser Verlag, Basel, (1989).

\bibitem{spr2} 
T.A. ~Springer, \emph{Linear Algebraic Groups},
Second edition. Progress in Mathematics, 9. Birkh\"auser Boston, Inc., 
Boston, MA, 1998.

\bibitem{St} 
R.~Steinberg, 
\emph{Endomorphisms of Linear Algebraic Groups},
Mem. Amer. Math. Soc. \textbf{80}, (1968).

\bibitem{St1} 
\bysame, 
\emph{Notes on Chevalley Groups}, Yale University, New Haven (1968).

\bibitem{SS}
T.A.~Springer, R.~Steinberg,
\emph{Conjugacy classes}, 
Seminar on algebraic groups and related finite groups, 
Lecture Notes in Mathematics, \textbf{131},  
Springer-Verlag, Heidelberg (1970), 167--266.

\bibitem{testerman}
D.M.~Testerman, 
\emph{A construction of certain maximal subgroups 
of the algebraic groups $E\sb 6$ and $F\sb 4$},
J. Algebra \textbf{122} (1989), no. 2, 299--322.

\end{thebibliography}
\end{document}